\documentclass{article}
\usepackage{amsfonts}
\usepackage{amsmath}
\usepackage{dsfont}
\usepackage{graphicx}
\usepackage{amssymb}
\usepackage{caption}
\usepackage[lined,boxed,commentsnumbered]{algorithm2e}
\usepackage{geometry}
\usepackage{csquotes}
\RequirePackage[OT1]{fontenc}
\newtheorem{theorem}{Theorem}

\newtheorem{remark}[theorem]{Remark}

\begin{document}

\title{Importance Sampling for multi-constraints rare event probability.}
\author{Virgile Caron\\Telecom ParisTech, 37-39 rue Dareau, 75014 Paris\\ virgile.caron@telecom-paristech.fr}
\maketitle

\begin{abstract}
Improving Importance Sampling estimators for rare event probabilities requires sharp approximations of the optimal density leading to a nearly zero-variance estimator. This paper presents a new way to handle the estimation of the probability of a rare event defined as a finite intersection of subset. We provide a sharp approximation of the density of long runs of a random walk conditioned by multiples constraints, each of them defined by an average of a function of its summands as their number tends to infinity.
\end{abstract}

\section{Introduction and context} \label{sec:intro}

In this paper, we consider efficient estimation of the probability of large deviations of a multivariate sum of independent, identically distributed, light-tailed and non-lattice random vectors. The probability to be estimated is defined as an intersection of $s$ event where each of them is characterized by a sum of $n$ i.i.d. random variable belonging to some countable union of intervals.  

Consider $\mathbf{X}_{1}^{n}:=\left(\mathbf{X}_{1},...,\mathbf{X}_{n}\right)$ $n$ i.i.d. random vectors with known common density $p_{\mathbf{X}}$ on $\mathbb{R}^{d}$, $d\geq{1},$ copies of $\mathbf{X}:=\left(\mathbf{X}^{(1)},...,\mathbf{X}^{(d)}\right).$ The superscript $(j)$ pertains to the coordinate of a vector and the subscript $i$ pertains to replications. Consider also $u$ a measurable function defined from $\mathbb{R}^{d}$ to $\mathbb{R}^{s}$ with $s<n$. The last condition on $s$ forbid to have a zero-one probability event. Define $\mathbf{U}:=u(\mathbf{X})$ with density $p_{\mathbf{U}}$ and 
\begin{equation*}
\mathbf{U}_{1,n}:=\sum_{i=1}^{n}\mathbf{U}_{i}.
\end{equation*}
We intend to estimate for large but fixed $n$ 
\begin{equation}\label{def:event_to_estime_simple}
P_{n}:=P\left( \mathbf{U}_{1,n}\in nA\right)
\end{equation}
where $A$ is a non-empty measurable set of $\mathbb{R}^{s}$ such as $E[u\left(\mathbf{X}\right)]\notin{A}.$  

For sake of clarity, the probability to be estimate can be written as follows
\begin{eqnarray}\label{def:event_to_estime}
P_{n}=P\left[\bigcap_{j=1}^{s}\left\{\sum_{i=1}^{n}u^{(j)}\left(\mathbf{X}_{i}^{(1)},...,\mathbf{X}_{i}^{(d)}\right)\in{nA^{(j)}}\right\}\right]
\end{eqnarray}
where, for $j\in\{1,...s\}$, $u^{(j)}$ is a measurable function from $\mathbb{R}^{d}$ to $\mathbb{R}$ and $A^{(j)}$ is a countable union of intervals of $\mathbb{R}.$ According to the context, we will use either (\ref{def:event_to_estime_simple}) or (\ref{def:event_to_estime}).
 
In \cite{BroniatowskiCaron2013}, the authors consider in details the case where $d=s=1$ and $A:=A_{n}=(a_{n},\infty )$ with $a_{n}$ a convergent sequence.

The basic estimate of $P_{n}$ is defined as follows: generate $L$ i.i.d.
samples $X_{1}^{n}(l)$ with underlying density $p_{\mathbf{X}}$ and define 
\begin{equation*}
\widetilde{P_{n}}:=\frac{1}{L}\sum_{l=1}^{L}\mathds{1}_{\mathcal{E}_{n}}\left(
X_{1}^{n}(l)\right)
\end{equation*}
where 
\begin{equation} \label{E_n}
\mathcal{E}_{n}:=\left\{ (x_{1},...,x_{n})\in \left(\mathbb{R}^{d}\right)^{n}:\left(
u\left(x_{1}\right)+...+u\left(x_{n}\right)\right)\in nA\right\}. 
\end{equation}
The Importance Sampling estimator of $P_{n}$ with sampling density $g$ on $\left(\mathbb{R}^{d}\right)^{n}$ is

\begin{equation}
\widehat{P_{n}}:=\frac{1}{L}\sum_{l=1}^{L}\hat{P}_{n}(l)\mathds{1}_{\mathcal{E%
}_{n}}\left( Y_{1}^{n}(l)\right)  \label{FORM_IS}
\end{equation}%
where $\hat{P}_{n}(l)$ is called "importance factor" and can be written 
\begin{equation} \label{Facteur_d_Importance}
\hat{P}_{n}(l):=\frac{\prod\limits_{i=1}^{n}p_{\mathbf{X}}\left(
Y_{i}(l)\right) }{g\left( Y_{1}^{n}(l)\right) } 
\end{equation}%
and where the $L$ samples $Y_{1}^{n}(l):=\left( Y_{1}(l),...,Y_{n}(l)\right)$ are i.i.d. with common density $g$; the coordinates of $Y_{1}^{n}(l)$ however need not be i.i.d.. It is known that the optimal choice for $g$ is the density of $\mathbf{X}_{1}^{n}:=\left( \mathbf{X}_{1},...,\mathbf{X}_{n}\right) $ conditioned upon $\left( \mathbf{X}_{1}^{n}\in \mathcal{E}_{n}\right) $, leading to a zero variance estimator. We refer to \cite{Bucklew2004} for the background of this section.

The state-independent IS scheme for rare event estimation (see \cite{BucklewNeySadowsky1990} or \cite{Sadowsky1996}), rests on two basic ingredients: the sampling distribution is fitted to the so-called dominating point (which is the point where the quantity to be estimated is mostly captured; see \cite{Ney1983}) of the set to be measured; independent and identically distributed replications under this sampling distribution are performed. More recently, a state-dependent algorithm leading to a strongly efficient estimator is provided by  \cite{Blanchetetal2008} when $d=1$, $u(x)=x$ and $A=(a;\infty)$ (or, more generally in $\mathbb{R}^{d}$, with a smooth boundary and a unique dominating point). Indeed, adaptive tilting defines a sampling density for the $i-$th r.v. in the run which depends both on the target event $\left( \mathbf{U}_{1,n}\in nA\right) $ and on the current state of the path up to step $i-1.$ Jointly with an ad hoc stopping rule controlling the excursion of the current state of the path, this algorithm provides an estimate of $P_{n}$ with a coefficient of variation independent upon $n$. This result shows that nearly optimal estimators can be obtained without approximating the conditional density.

The main issue of the method described above is to find dominating point. However, when the dimension of the set $A$ increases, finding a dominating point can be very tricky or even impossible. A solution will be to divide the set under consideration into smaller subset and, for each one of this subset, find a dominating point. Doing so makes the implementation of an IS scheme harder and harder as the dimension increases.

Our proposal is somehow different since it is based on a sharp approximation
result of the conditional density of long runs. The approximation holds for
any point conditioning of the form $\left( \mathbf{U}_{1,n}=nv\right).$ Then
sampling $v$ in $A$ according to the distribution of $\mathbf{U}_{1,n}$
conditioned upon $\left( \mathbf{U}_{1,n}\in nA\right) $ produces the
estimator. By its very definition this procedure does not make use of any
dominating point, since it randomly explores the set $A.$ Indeed, our proposal hints on two choices: first do not make use of the notion of dominating point and explore all the target set instead (no part of the set $A$ is neglected); secondly, do not use i.i.d. replications, but merely sample long runs of variables under a proxy of the optimal sampling scheme. 

We will propose an IS sampling density which approximates this conditional
density very sharply on its first components $y_{1},...,y_{k}$ where $k=k_{n} $ is very large, namely $k/n\rightarrow 1.$ However, but in the Gaussian case, $k$ should satisfy $\left(n-k\right)\rightarrow\infty$ by the very construction of the approximation. The IS density on $\left(\mathbb{R}^{d}\right)^{n}$ is obtained multiplying this proxy by a product of a much simpler state-independent IS scheme following \cite{SadowskyBucklew1990}.

The paper is organized as follows. Section \ref{sec:notations_hypothesis} is devoted to notations and hypothesis. In Section \ref{sec:approx_local}, we expose the approximation scheme for the conditional density of $\mathbf{X}_{1}^{k}$ under $\left(\mathbf{U}_{1,n}=nv\right).$  In Section \ref{sec:define_estimator}, we introduce our IS scheme. Simulated results are presented in Section \ref{sec:example} which enlighten the gain of the present approach over state-dependent Importance Sampling schemes.

We rely on \cite{Caron2013} where the basic approximation (and proofs) used
in the present paper can be found. The real case is studied in \cite{BroniatowskiCaron2014} and applications for IS estimators can be found in \cite{BroniatowskiCaron2013}.

\section{Notations and hypotheses} \label{sec:notations_hypothesis}

We consider approximations of the density of the vector $\mathbf{X}_{1}^{k}$ on $\left(\mathbb{R}^{d}\right)^{k}$, when the conditioning event writes (\ref{def:event_to_estime}) and $k:=k_{n}$ is such that 
\begin{equation} \label{cond:K1}
\lim_{n \to \infty} \frac{k}{n}=1  \tag{K1}
\end{equation}
\begin{equation} \label{cond:K2}
\lim_{n \to \infty} (n-k)=+\infty. \tag{K2}
\end{equation}
Therefore we may consider the asymptotic behavior of the density of the random walk on long runs.

Throughout the paper the value of a density $p_{\mathbf{Z}}$ of some
continuous random vector $\mathbf{Z}$ at point $z$ may be written $p_{\mathbf{Z}}(z)$ or $p\left( \mathbf{Z}=z\right) ,$ which may prove more
convenient according to the context. 

Let $p_{nv}$ denote the density of $\mathbf{X}_{1}^{k}$ under the local
condition $\left(\mathbf{U}_{1,n}=nv\right)$ 
\begin{equation}  \label{def:p_cond_ponc_chap2}
p_{nv}\left( \mathbf{X}_{1}^{k}=Y_{1}^{k}\right) :=p(\left. \mathbf{X}%
_{1}^{k}=Y_{1}^{k}\right\vert \mathbf{U}_{1,n}=nv)
\end{equation}
where $Y_{1}^{k}$ belongs to $\left(\mathbb{R}^{d}\right)^{k}$ and $v$ fixed in $\mathbb{R}^{s}.$

We will also consider the density $p_{nA}$ of $\mathbf{X}_{1}^{k}$
conditioned upon $\left( \mathbf{U}_{1,n}\in nA\right) $ 
\begin{equation}  \label{def:p_cond_plusgrand_chap2}
p_{nA}\left( \mathbf{X}_{1}^{k}=Y_{1}^{k}\right) :=p(\left. \mathbf{X}%
_{1}^{k}=Y_{1}^{k}\right\vert\mathbf{U}_{1,n}\in nA).
\end{equation}

The approximating density of $p_{nv}$ is denoted $g_{nv}$; the corresponding
approximation of $p_{nA}$ is denoted $g_{nA}.$ Explicit formulas for those
densities are presented in the next section.

\section{Multivariate random walk under a local conditionning event.} \label{sec:approx_local}

Let $\epsilon_{n}$ be a postive sequence such as 
\begin{equation} \label{cond:A1}
\lim_{n \to \infty} \epsilon_{n}^{2}(n-k)=\infty \tag{E1}
\end{equation}
\begin{equation} \label{cond:A2}
\lim_{n \to \infty} \epsilon_{n}(\log n)^{2}=0 \tag{E2}
\end{equation}

It will be shown that $\epsilon_{n}\left( \log n\right) ^{2}$ is the rate of
accuracy of the approximating scheme.

We assume that $\mathbf{U}:=u\left(\mathbf{X}\right)$ has a density $p_{\mathbf{U}}$ (with p.m. $P_{\mathbf{U}}$) absolutely continuous with respect to Lebesgue measure on $\mathbb{R}^{s}.$ Futhermore, we assume that $u$ is such that the characteristic function of $\mathbf{U}$ belongs to $L^{r}$ for some $r\geq{1}.$

Denote $\underline{0}$ is the vector of $\mathbb{R}^{s}$ with all coordinates equal to $0$ and $V(\underline{0})$ a neighborhood of $\underline{0}.$

We assume that $\mathbf{U}$ satisfy the Cramer condition, meaning
\[\Phi_{\mathbf{U}}(t):=E[\exp<t,\mathbf{U}>]<\infty,\textsc{\ \ } t\in V(\underline{0})\subset\mathbb{R}^{s}.\]
and denote 
\begin{eqnarray*}
m(t):= {}^{t}\nabla \log(\Phi_{\mathbf{U}}(t)),\textsc{\ \ } t\in V(\underline{0})\subset\mathbb{R}^{s}
\end{eqnarray*}
and
\begin{eqnarray*}
\kappa(t):= {}^{t}\nabla \nabla \log(\Phi_{\mathbf{U}}(t)),\textsc{\ \ } t\in V(\underline{0})\subset\mathbb{R}^{s}.
\end{eqnarray*}
as the mean and the covariance matrix of the tilted density defined by
\begin{eqnarray} \label{def:pi_non_centre_r_d}
\pi_{u}^{\alpha}(x) :=\frac{\exp<t,u(x)>}{\Phi_{\mathbf{U}}(t)}p_{\mathbf{X}}(x).
\end{eqnarray}
where $t$ is the only solution of $m(t)=\alpha$ for $\alpha$ in the conxev hull of $P_{\mathbf{U}}.$ Conditions on $\Phi_{\mathbf{U}}(t)$ which ensure existence and uniqueness of $t$ are referred to stepness properties (see \cite{BarndorffNielsen1978}, p153 ff, for all properties of moment generating funtion used in this paper).

Let $v\in\mathbb{R}^{s}.$ We now state the general form of the approximating density.
Denote
\begin{eqnarray}
g_{0}(y_{1}|y_{0}):=\pi_{u}^{v}(y_{1})
\end{eqnarray}
with an arbitrary $y_{0}$ and $\pi_{u}^{v}$ defined in (\ref{def:pi_non_centre_r_d}).

For $1\leq{i\leq{k-1}}$, we recursively define $g(y_{i+1}|y_{1}^{i})$. Set $t_{i}\in\mathbb{R}^{s}$ to be the unique solution to the equation
\begin{eqnarray}
m(t_{i})=m_{i,n}:=\frac{n}{n-i}\left(v-\frac{u_{1,i}}{n}\right)
\end{eqnarray}
where $u_{1,i}=u(y_{1})+...+u(y_{i}).$

Denote
\begin{eqnarray*}
\kappa_{(i,n)}^{j,l}:=\frac{d^{2}}{dt^{(j)}dt^{(l)}}\left(\log E_{\pi_{\mathbf{U}}^{m_{i,n}}}\exp <t,\mathbf{U}>\right)\left(0\right) 
\end{eqnarray*}
and 
\begin{eqnarray*}
\kappa_{(i,n)}^{j,l,m}:=\frac{d^{3}}{dt^{(j)}dt^{(l)}dt^{(m)}}\left(\log E_{\pi_{\mathbf{U}}^{m_{i,n}}}\exp <t,\mathbf{U}>\right)\left(0\right) .
\end{eqnarray*}
for $j,l$ and $m$ in $\{1,...,s\}$

Denote
\begin{eqnarray}
g(y_{i+1}|y_{1}^{i}):=C_{i}\mathfrak{n}_{s}\left(u(y_{i+1});\beta\alpha+v,\beta\right)p_{\mathbf{X}}(y_{i+1})
\end{eqnarray}
where $C_{i}$ is a normalizing factor, $\mathfrak{n}_{s}\left(u(y_{i+1});\beta\alpha+v,\beta\right)$ is the normal density at $u(y_{i+1})$ with mean $\beta\alpha+v$ and covariance matrix $\beta$ with $\alpha$ and $\beta$ defined by
\[
\alpha:=\left(t_{i}+\frac{\kappa_{(i,n)}^{-2}\gamma}{2(n-i-1)}\right)\\
\beta:=\kappa_{(i,n)}(n-i-1)
\]
where $\gamma$ is defined by 
\[
\gamma:=\left(\sum_{j=1}^{s} \kappa_{(i,n)}^{j,j,p}\right)_{1\leq{p}\leq{s}}.
\]

Then
\begin{eqnarray}\label{def:g_nv}
g_{nv}(y_{1}^{k}):=g_{0}(y_{1}|y_{0})\prod_{i=1}^{k-1} g(y_{i+1}|y_{1}^{i})
\end{eqnarray}

\begin{theorem} \label{th:princ_general}
Assume (\ref{cond:A1}), (\ref{cond:A2}), (\ref{cond:K1}) and (\ref{cond:K2}).
\begin{itemize}
\item Let $Y_{1}^{k}$ a sample of $P_{nv}.$
Then
\begin{eqnarray} \label{th:princ_general_1}
p\left(\textbf{X}_{1}^{k}=Y_{1}^{k}|\mathbf{U}_{1,n}=nv\right)=g_{nv}(Y_{1}^{k})
(1+o_{P_{nv}}(1+\epsilon_{n}(\log n)^{2}))
\end{eqnarray}
\item Let $Y_{1}^{k}$ a sample of $G_{nv}.$
Then
\begin{eqnarray} \label{th:princ_general_inverse}
p\left(\textbf{X}_{1}^{k}=Y_{1}^{k}|\mathbf{U}_{1,n}=nv\right)=g_{nv}(Y_{1}^{k})
(1+o_{G_{nv}}(1+\epsilon_{n}(\log n)^{2}))
\end{eqnarray}
\end{itemize}
\end{theorem}

\begin{remark}
The approximation of the density of $\mathbf{X}_{1}^{k}$ is not performed on the sequence of entire spaces $\left(\mathbb{R}^{s}\right)^{k}$ but merely on a sequence of subsets of $\left(\mathbb{R}^{s}\right)^{k}$ which contains the trajectories of the conditioned random walk with probability going to $1$ as $n$ tends to infinity. The approximation is performed on \textit{typical paths}. For sake of applications in Importance Sampling, Theorem \ref{th:princ_general_inverse} is exactly what we need. Nevertheless, as proved in \cite{Caron2013}, the extension of our results from typical paths to the whole space $\mathbb{R}^{k}$ holds: convergence of the relative error on large sets imply that the total variation distance between the conditioned measure and its approximation goes to $0$ on the entire space. 
\end{remark}

\begin{remark}
The rule which defines the value of $k$ for a given accuracy of the approximation is stated in Section 5 of \cite{Caron2013}.
\end{remark}

\begin{remark} \label{remark:gaussian_exact}
When the $\mathbf{X}_{i}$'s are i.i.d. Gaussian standard and $u(x)=x$, the result of the approximation theorem are true for $k=n-1$ without the error term. Indeed, it holds $p(\left. \mathbf{X}_{1}^{n-1}=x_{1}^{n-1}\right\vert \mathbf{U}_{1,n}=nv)=g_{nv}\left( x_{1}^{n-1}\right) $ for all $x_{1}^{n-1}$ in $\left(\mathbb{R}^{d}\right)^{n-1}$.
\end{remark}

As stated above the optimal choice for the sampling density is $p_{nA}.$ It holds 
\begin{equation}
p_{nA}(x_{1}^{k})=\int_{A}p_{nv}\left( \mathbf{X}%
_{1}^{k}=x_{1}^{k}\right) p(\left. \mathbf{U}_{1,n}/n=v\right\vert \mathbf{U}%
_{1,n}\in nA)dv  \label{etoile_chap2}
\end{equation}%
so that, in contrast with \cite{Blanchetetal2008} or \cite{BucklewNeySadowsky1990}, we do not consider the dominating point approach but merely realize a sharp
approximation of the integrand at any point of $A$ and consider the dominating contribution of all those distributions in the evaluation of the conditional density $p_{nA}.$

\section{Adaptive IS Estimator for rare event probability} \label{sec:define_estimator}

The IS scheme produces samples $Y:=\left( Y_{1},...,Y_{k}\right) $
distributed under $g_{nA}$, which is a continuous mixture of densities $g_{nv}$ as in (\ref{def:g_nv}) with $p\left(\mathbf{U}_{1,n}/n=v\vert\mathbf{U}_{1,n}\in nA\right).$

Simulation of samples $\mathbf{U}_{1,n}/n$ under this density can be
performed through Metropolis-Hastings algorithm, since 
\begin{equation*}
r(v,v^{\prime }):=\frac{p(\left. \mathbf{U}_{1,n}/n=v\right\vert \mathbf{U}%
_{1,n}\in nA)}{p(\left. \mathbf{U}_{1,n}/n=v^{\prime }\right\vert \mathbf{U}%
_{1,n}\in nA)}
\end{equation*}%
turns out to be independent upon $P\left( \mathbf{U}_{1,n}\in nA\right) .$
The proposal distribution of the algorithm should be supported by $A.$

The density $g_{nA}$ is extended from $\left(\mathbb{R}^{s}\right)^{k}$ onto $\left(\mathbb{R}^{s}\right)^{n}$
completing the $n-k$ remaining coordinates with i.i.d. copies of r.v's $%
Y_{k+1},...,Y_{n}$ with common tilted density 
\begin{equation}
g_{nA}\left( \left. y_{k+1}^{n}\right\vert y_{1}^{k}\right) :=\prod
\limits_{i=k+1}^{n}\pi_{u}^{m_{k}}(y_{i})  \label{complementk,n}
\end{equation}
with $m_{k}:=m(t_{k})=\frac{n}{n-k}\left( v-\frac{u_{1,k}}{n}\right)$ and
\begin{equation*}
u_{1,k}=\sum_{i=1}^{k}u(y_{i}).
\end{equation*}

The last $n-k$ r.v's $\mathbf{Y}_{i}$'s are therefore drawn according to the
state independent i.i.d. scheme in phase with Sadowsky and Bucklew \cite%
{SadowskyBucklew1990}.

We now define our IS estimator of $P_{n}.$ Let $Y_{1}^{n}(l):=Y_{1}(l),...,Y_{n}(l)$ be generated under $g_{nA}.$ Let 
\begin{equation}
\widehat{P_{n}}(l):=\frac{\prod_{i=0}^{n}p_{\mathbf{X}}(Y_{i}(l))}{%
g_{nA}(Y_{1}^{n}(l))}\mathds{1}_{\mathcal{E}_{n}}\left( Y_{1}^{n}(l)\right)
\label{Pchapl}
\end{equation}%
and define 
\begin{equation}
\widehat{P_{n}}:=\frac{1}{L}\sum_{l=1}^{L}\widehat{P_{n}}(l).  \label{Pchap}
\end{equation}%
in accordance with (\ref{FORM_IS}).

\begin{remark} 
In the real case and for $A=(a,\infty)$, the authors of \cite{BroniatowskiCaron2013} shows that under regularity conditions the resulting relative error of the estimator is proportional to $\sqrt{n-k_{n}}$ and drops by a factor $\sqrt{n-k_{n}}/\sqrt{n}$ with respect to the state independent IS scheme. Slight modification in the extension of $g_{nA}$ allows to prove the strong efficiency of the estimator (\ref{Pchap}) using arguments from both \cite{Blanchetetal2008} and \cite{BroniatowskiCaron2013}; see \cite{Caronetal2013}.
\end{remark}

\section{When the dimension becomes very high} \label{sec:example}

This section compares the performance of the present approach with respect
to the standard tilted one using i.i.d. replications under (\ref{def:pi_non_centre_r_d}). Consider a random sample $X_{1},...,X_{100}$ where $X_{1}$ has a normal distribution $N(0.05,1)$ and let 
\begin{equation*}
\mathcal{E}_{100}:=\left\{ x_{1}^{100}:\frac{\left\vert
x_{1}+...+x_{100}\right\vert }{100}>0.28\right\}.
\end{equation*}

This example is in the same vein as the one developed in \cite{DupuisWang2004} or in \cite{GlassermanWang1997}. Under the present proposal the distribution of the Importance Factor concentrates around $P_{100}$; hence so-called "rogue path phenomenon" (see \cite{DupuisWang2004}) does not occur.

We explore the gain in relative accuracy when the dimension of the measured set increases. Let therefore $B:=\left( \mathcal{E}_{100}\right)^{d}$ which is the $d$-cartesian product of $\mathcal{E}_{100}.$ The $100$ r.v.'s $X_{i}$ 's are i.i.d. random vectors in $\mathbb{R}^{d}$ with common i.i.d. $N(0.05,1)$ distribution. The dominating point has all coordinates equal 0.28. Rogue path curse produces an overwhelming loss in accuracy, imposing a very large increase in runtime to get reasonable results. Our interest is to show that in this simple dissymetric case our proposal provides a good estimate, while the standard IS scheme ignores a part of the event $B.$ The standard i.i.d. IS scheme introduces the dominating point $a=0.28$ and the family of i.i.d. tilted r.v's with common $N(a,1)$ distribution for each coordinates. It can be seen that a large part of $B$ is never visited through the procedure, inducing a bias in the estimation.

This example is not as artificial as it may seem; indeed it leads to a $2^{d}$
dominating points situation which is quite often met in real life. Exploring
at random the set of interest avoids any search for dominating points. Drawing $L$ i.i.d. points $v_{1},...,v_{L}$ according to the distribution of $\mathbf{U}_{1,100}/100$ conditionally upon $B$ we evaluate $P_{100}$ with $k=99$; note that in the Gaussian case Theorem \ref{th:princ_general} provides an exact description of the conditional density of $X_{1}^{k}$ for all $k$ between $1$ and $n$. The following figure shows the gain in relative accuracy w.r.t. the state independent IS scheme according to the growth of $d.$ The value of $P_{100}$ is $10^{-2d}.$

\begin{figure}[th]
\centering \includegraphics*[scale=0.35]{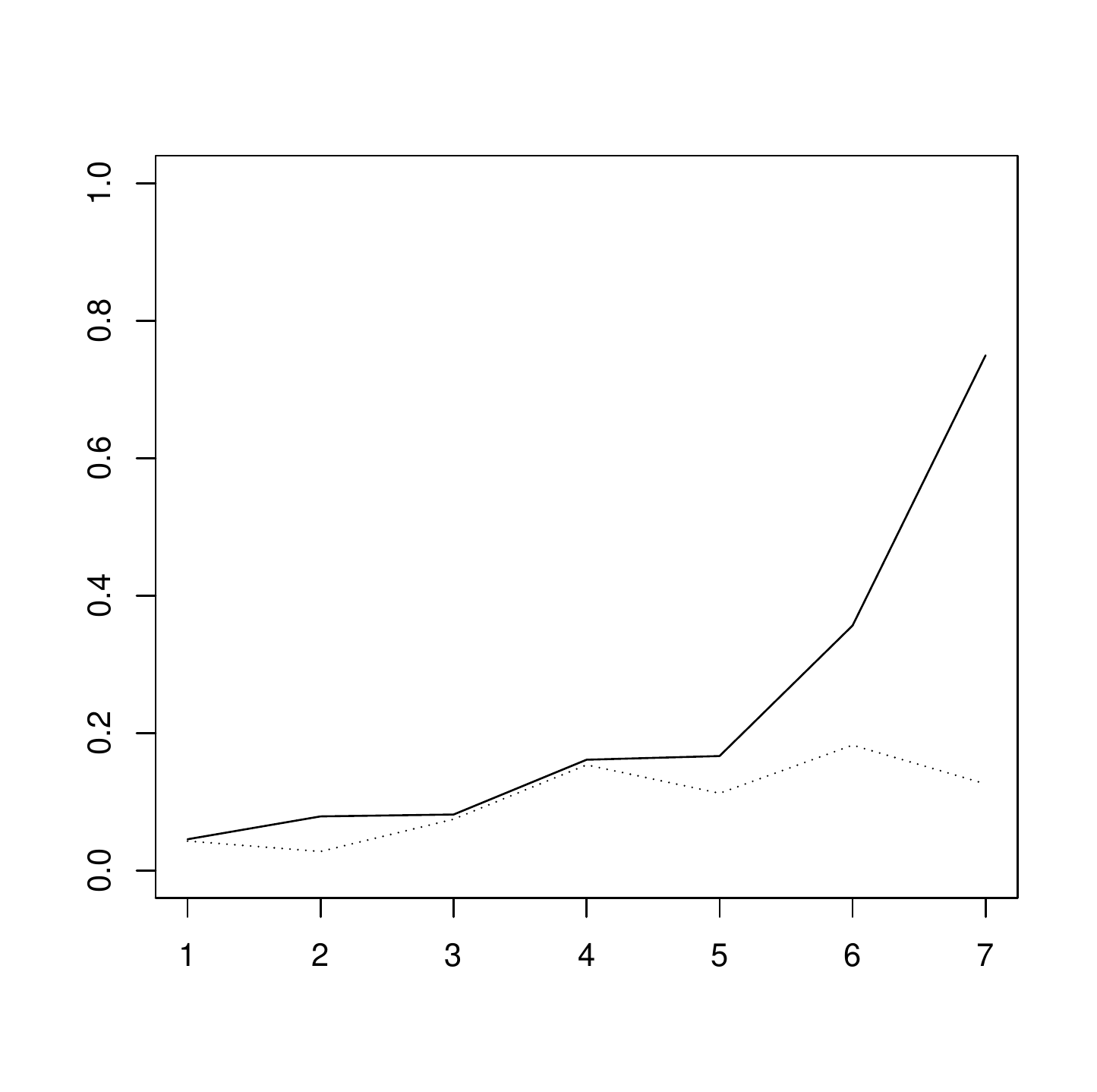}
\caption{Relative Accuracy of the adaptive estimate (dotted line) w.r.t.
i.i.d. twisted one (solid line) as a function of the dimension $d$ for $L=1000.$}
\end{figure}

\section{Conclusion}

In this paper, we explore a new way to estimate multi-constraints large deviation probability. In future work, the author will investigate the theoretical behaviour of the relative error of our proposed estimator.

\end{document}